\documentclass[12pt]{amsart}
\usepackage{graphicx}
\usepackage{amsmath}
\usepackage{latexsym}
\usepackage{amsthm}
\usepackage{amsfonts}
\usepackage{amssymb}
\usepackage[all]{xy}

\theoremstyle{plain} 
\newtheorem{theorem}{Theorem}[section]
\newtheorem{lemma}[theorem]{Lemma}
\newtheorem{proposition}[theorem]{Proposition}
\newtheorem{corollary}{Corollary}[theorem]

\begin{document}
\title{Recursion categories of coalgebras}
\author{Florian Lengyel}
\address{CUNY Graduate Center, 365 5th Avenue, New York, NY 10016}
\email{flengyel@gc.cuny.edu}
\begin{center}
\thanks{Typeset with \LaTeX$2\epsilon$ and \Xy-pic version \xyversion}
\end{center}
\begin{abstract}
We construct recursion categories from categories of coalgebras.
Let $F$ be a nontrivial endofunctor on the category of sets that
weakly preserves pullbacks and such that the category
$\textbf{Set}_F$ of $F$-coalgebras is complete. The category
$\textbf{Set}_F$ may be embedded in the category $\mathbf{Pfn}_F$
of $F$-coalgebras and partial morphisms, which is a $P$-category
that is prodominical but not dominical in general. An existence
theorem of A. Heller is applied to certain subcategories of
$\textbf{Pfn}_F$ to obtain examples of recursion categories of
coalgebras.
\end{abstract}

\maketitle

\section{Introduction}
The theory of universal coalgebras has recently undergone vigorous
development as a descriptive study of automata, data structures,
and the semantics of object oriented programming, among other
applications ~\cite{MR1791953,glutacs}. Coalgebras as data
structures include colorings, permutations, trees, partial
self-mappings, the real numbers, streams and formal power series
~\cite{glutacs,MR2000d:03174,MR1659609}. Dominical categories
~\cite{MR88i:03078} were developed by DiPaola and Heller  as an
element-free generalization of classical recursion theory, in
which one has categorical analogs of the construction of the
universal Turing machine, the s-m-n theorem and the Kleene
recursion theorem. Rosolini's $P$-categories ~\cite{Rosolini}
generalize the dominical categories and have gained acceptance as
an appropriate setting for categorical developments of recursion
theory and programming semantics
~\cite{MR92e:03067,TaylorP:exaiw}. We construct prodominical
recursion categories of coalgebras; in particular, this shows that
many categories of coalgebras possess a Turing morphism, which is
a categorical generalization of the universal Turing machine.

Let $F$ be a nontrivial endofunctor on the category of sets such
that $F$ weakly preserves pullbacks, and such that the category
$\textbf{Set}_F$ of $F$-coalgebras is complete. For example, if
the forgetful functor $\textbf{Set}_F\rightarrow\textbf{Set}$ has
a right adjoint, then $\textbf{Set}_F$ is complete ~\cite{Kurz}.
Under various mild assumptions on the endofunctor $F$, such as
boundedness or accessiblility, the right adjoint will exist
~\cite{GS2000,Kurz,MR99h:18012}. The category $\textbf{Set}_F$ may
be embedded in the category $\mathbf{Pfn}_F$ of $F$-coalgebras and
partial morphisms, which is a $P$-category that is prodominical
but not dominical in general.

This paper is organized as follows. In the first section we recall
definitions and results needed from the theory of coalgebras that
we use to construct the category $\mathbf{Pfn}_F$ of coalgebras
and partial morphisms. Next, we review A. Heller's theory of
recursion categories ~\cite{MR92e:03067}, and verify that
$\mathbf{Pfn}_F$ satisfies certain properties required for the
application of Heller's existence theorem. In the second section
we state Heller's existence theorem and apply it to obtain our
existence result.  In the final section we show that
$\mathbf{Pfn}_F$ contains many recursion categories. In
conclusion, we mention related results beyond the scope of this
paper.

\subsection{Coalgebras} Let $F$ be an endofunctor on the
category $\mathbf{Set}$ of sets. An {\it $F$-coalgebra} $(A,
\alpha)$ is a set $A$ together with a map $\alpha:A\rightarrow F
A$ called its {\it structure map}. A morphism of $F$-coalgebras
$\varphi:(A,\alpha)\rightarrow(B,\beta)$ is a map
$\varphi:A\rightarrow B$ such that
$\beta\circ\varphi=F\varphi\circ\alpha$. The $F$-coalgebras and
their morphisms form a category, denoted by $\mathbf{Set}_F$.
Coalgebras in this sense are more general than those associated
with comonads ~\cite{MR1712872}. The category ${\mathbf Set}_F$ is
cocomplete, co-wellpowered, and is closed under homomorphic images
~\cite{glutacs}.

Let $(A,\alpha)$ be a coalgebra. A subset $U$ of $A$ has at most
one coalgebra structure such that the inclusion $U\hookrightarrow
A$ induces a morphism of coalgebras. A {\it subcoalgebra} of a
coalgebra $(A,\alpha)$ is a coalgebra $(U,\beta)$ such that the
inclusion is a morphism. The set endofunctor $F$ is {\it
nontrivial} if for any set $X$, $FX=\emptyset$ implies
$X=\emptyset$. If $F$ is nontrivial, then the set of subcoalgebras
of a given $F$-coalgebra $X$ forms a topology on $X$ called the
{\it coalgebra topology} ~\cite{GS2000}.

A commutative square in a category $\mathbf{C}$ is a {\it weak
pullback} if it satisfies all of the conditions of a pullback
except for the uniqueness requirement for induced morphisms.
Analogously, one may define a {\it weak limit} of a diagram in
$\mathbf{C}$. A {\it sink} is a collection of morphisms of
$\mathbf{C}$ with a common codomain. A {\it (weak) generalized
pullback} is a (weak) limit of a sink. A functor {\it preserves
weak (generalized) pullbacks} if it sends weak (generalized)
pullbacks to weak (generalized) pullbacks. A functor {\it weakly
preserves (generalized) pullbacks} if it sends (generalized)
pullbacks to weak (generalized) pullbacks. If (generalized)
pullbacks exist in $\mathbf{C}$, and if
$F:\mathbf{C}\rightarrow\mathbf{Set}$ is a set valued functor,
then $F$ preserves weak (generalized) pullbacks if and only if $F$
weakly preserves (generalized) pullbacks ~\cite{glutacs}. More
generally, if $\mathbf{C}$ has all limits of certain diagrams $D$,
then $F$ preserves weak limits of $D$-diagrams if and only if $F$
weakly preserves limits of $D$-diagrams.

The following $\mathbf{Set}$ endofunctors preserve weak limits of
sinks: the identity functor, all constant functors, functors of
the form $X\mapsto X\coprod Y$ for a fixed set $Y$, functors of
the form $X\mapsto X^Y$ for a fixed set $Y$, and the powerset
functor. Moreover, sums, products and composites of $\mathbf{Set}$
endofunctors that preserve weak limits of sinks also preserve weak
limits of sinks ~\cite{MR1787576}.

Gumm  and Schr\"oder  ~\cite{MR1787576} prove that the
$\mathbf{Set}$ endofunctor $F$ weakly preserves pullbacks along
monomorphisms if and only for any morphism $\varphi:A\rightarrow
B$ of $F$-coalgebras and for any subcoalgebra $U$ of $B$, the
preimage $\varphi^{-1}[U]$ of $U$ under $\varphi$ is a
subcoalgebra of $A$.

A $\mathbf{Set}$ endofunctor $F$ is {\it bounded by the cardinal
$\kappa$} if for each $F$-coalgebra $A$, and for each $a\in A$
there is a subcoalgebra $B$ of $A$ containing $a$ such that
$|B|\leq\kappa$. Let $C$ and $M$ be two fixed, possibly empty
sets. Any functor of the form $X\mapsto C\times X^M$ is bounded;
in particular, the identity functor and the constant functors are
bounded. Moreover, sums, products and compositions of bounded
functors are bounded. A characterization of bounded functors is
given in ~\cite{GSbounded}. Examples of coalgebras for bounded
functors include colorings, which are coalgebras for a constant
$\mathbf{Set}$ endofunctor, self-maps of a set, which are
coalgebras for the identity functor, and stream coalgebras, which
are coalgebras for the $\mathbf{Set}$ endofunctor given by
$X\mapsto $A$\times X$, where $A$ is a fixed nonempty set
~\cite{glutacs,MR1659609,MR2000d:03174}. In each of these
examples, the type endofunctor preserves weak limits of sinks.

Let $U:\mathbf{Set}_F\rightarrow\mathbf{Set}$ be the {\it
forgetful} functor, which sends an $F$-coalgebra
$\alpha:A\rightarrow F(A)$ to the set $A$, and a morphism of
coalgebras to itself as a morphism of $\mathbf{Set}$.  Kurz
~\cite{Kurz} has shown that if the forgetful functor $U$ has a
right adjoint, then $\mathbf{Set}_F$ is complete. The existence of
the right adjoint implies that the cofree coalgebra on any set
exists; limits are constructed as certain subcoalgebras of cofree
coalgebras. Gumm and Schr\"oder construct such a right adjoint to
the forgetful functor $U$ in the case when the endofunctor $F$ is
bounded; hence if $F$ is bounded, $\mathbf{Set}_F$ is complete
~\cite{GS2000}.

Let $X_i, i=0,1$ be $F$-coalgebras. A {\it bisimulation} between
$X_0$ and $X_1$ is a relation $R\subseteq U X_0\times U X_1$
equipped with an $F$-coalgebra structure $\rho:R\rightarrow F(R)$,
called a {\it bisimulation structure}, such that the projections
$\pi_i:R\rightarrow X_i,i=0,1$ are morphisms of coalgebras
~\cite{glutacs,MR1791953}. A bisimulation structure need not be
unique if it exists. If $X$ is a coalgebra, and if
$\varphi_i:X\rightarrow Y_i,i=0,1$ are coalgebra morphisms, then
$\left(U\varphi_0\times U\varphi_1\right)U X$ is a bisimulation
between $Y_0$ and $Y_1$ ~\cite{glutacs}. Let $X,Y$ be coalgebras
and let $\varphi:U X\rightarrow U Y$ be a map in $\mathbf{Set}$.
Then $\varphi$ defines a coalgebra morphism $X\rightarrow Y$ if
and only if its graph is a bisimulation; it follows that for any
coalgebra $X$, the diagonal $\Delta_X=\{(x,x):x\in U X\}$ is a
bisimulation ~\cite{MR1791953}.

A functor between $\kappa$-accessible categories is {\it
$\kappa$-accessible} if it preserves $\kappa$-filtered colimits.
It follows from results of Power and Watanabe ~\cite{MR99h:18012}
that if $F$ is an $\omega$-accessible endfunctor, then
$\mathbf{Set}_F$ is complete.

\subsection{Partial morphisms of coalgebras}
Let $F$ denote an endofunctor on $\mathbf{Set}$. Under certain
restrictions on $F$, the category $\mathbf{Set}_F$ can be extended
to the category of $F$-coalgebras and partial morphisms, denoted
by $\mathbf{Pfn}_F$, which we will obtain by following the
procedure in Rosolini's thesis ~\cite{Rosolini} for constructing a
category of partial maps.

Let $\mathbf{C}$ be a category with finite limits. A {\it
dominion} is a class $\mathcal{M}$ of monomorphisms of
$\mathbf{C}$ that is closed under identities, composition and
pullbacks of morphisms of $\mathbf{C}$ ~\cite{Rosolini}. If
$\mathcal{M}$ is a dominion of $\mathbf{C}$ and if $X$ and $Y$ are
objects of $\mathbf{C}$, a {\it partial map from $X$ to $Y$
defined in $\mathcal{M}$} is a pair $(m,\varphi)$, where
$m:U\rightarrow X$ is in $\mathcal{M}$, and where
$\varphi:U\rightarrow Y$ is a morphism of $\mathbf{C}$. Two such
pairs $(m,\varphi)$ and $(m',\varphi')$ are {\it equivalent} if
there exists an isomorphism $\theta$ making the following obvious
diagram commute.
\[
\newdir{ >}{{}*!/-10pt/@{>}}
\xymatrix{ U \ar[r]^-\varphi\ar@{ >->}[d]_-m\ar@{.>}[dr]^-\theta & Y \\ %
X & U' \ar@{ >->}[l]^-{m'} \ar[u]_-{\varphi'} }
\]
The equivalence class of the pair $(m,\varphi)$ is denoted by
$\{m,\varphi\}$ and is called a {\it partial morphism from $X$ to
$Y$}, also denoted by $X\dashrightarrow Y$. Composition of partial
morphisms $A\dashrightarrow B$ and $B\dashrightarrow C$ is defined
by the following diagram in which the square is a pullback.
\[
\newdir{ >}{{}*!/-10pt/@{>}}
\xymatrix{ W \ar[r]\ar@{ >->}[d] & V \ar[r]^{\varphi'}\ar@{ >->}[d]_{m'}& C \\%
U \ar[r]^\varphi\ar@{ >->}[d]_m & B\ar@{-->}[ur]_{\{m',\varphi'\}}\\ %
A\ar@{-->}[ur]_{\{m,\varphi\}} }
\]
The equivalence class $\{m'',\varphi''\}$ of the composite is
determined by the monomophism $m'':W\rightarrow A$ obtained as the
composite of the morphisms of the left vertical column, and the
morphism $\varphi'':W\rightarrow C$ obtained as the composite of
the morphisms of the top horizontal row. We can then express the
composite by the equation
\[
\{m'',\varphi''\}=\{m',\varphi'\}\circ\{m,\varphi\}.
\]
The composite is associative. For each object $A$, the equivalence
class $\{1_A,1_A\}$ acts as an identity for composition. Given a
category $\mathbf{C}$ with finite products together with a
dominion $\mathcal{M}$, the {\it category of partial maps
$P(\mathbf{C},\mathcal{M})$} is the category whose objects are
those of $\mathbf{C}$ and whose morphisms are the equivalence
classes of partial maps defined in $\mathcal{M}$.

We will define the category $\mathbf{Pfn}_F$ as
$P(\mathbf{Set}_F,\mathcal{M})$ by taking the dominion
$\mathcal{M}$ to be the class of monomorphisms of $\mathbf{Set}_F$
subject to conditions on the endofunctor $F$ that will make
$\mathbf{Set}_F$ complete and such that monomorphisms will
correspond to subcoalgebras.

For the first requirement, we note that if the forgetful functor
$\mathbf{Set}_F\rightarrow\mathbf{Set}$ has a right adjoint, which
happens, for example, if $F$ is bounded or accessible, then
$\mathbf{Set}_F$ is complete ~\cite{GS2000, Kurz, MR99h:18012}; in
particular, $\mathbf{Set}_F$ has products.

For the second requirement, we note that although a morphism in
$\mathbf{Set}_F$ that is injective in $\mathbf{Set}$ is a
monomorphism in $\mathbf{Set}_F$, the converse does not hold in
general ~\cite{glutacs,MR1791953}, and so monomorphisms do not
necessarily correspond to subcoalgebras in $\mathbf{Set}_F$. To
obtain a workable notion of a partial map, we will use the
following result on factorization systems in $\mathbf{Set}_F$ due
to Kurz, who introduced in his thesis the technique of lifting
factorization systems from the base category to obtain a
characterization of subcoalgebras in general categories of
coalgebras ~\cite{Kurz}.\footnote{It is sufficient for our
application to cite Proposition 4.7 of ~\cite{MR1791953}, which
states that if $F$ preserves weak pullbacks, then monomorphisms of
$\mathbf{Set}_F$ are precisely the injective morphisms; however,
the use of factorization systems applies to more general
categories of coalgebras.}

Let $\mathbf{C}$ be a category. A monomorphism is {\it regular}
if it is an equalizer. A monomorphism $m$ in $\mathbf{C}$ is {\it
extremal} if and only if $m=f e$ for some $f$ and for some epi
$e$ implies that $e$ is an isomorphism. A monomorphism $m$ is
{\it strong} if and only if for all epis $e$ and for all $f,g$
such that the following square commutes, there is a unique $d$
such that the triangles in the following diagram commute.
\[
\newdir{ >}{{}*!/-10pt/@{>}}
\xymatrix{ A \ar@{>>}[r]^e\ar[d]_f &B\ar@{.>}[dl]_d\ar[d]^g\\ %
C \ar@{ >->}[r]_m &D }
\]
Strong monomorphisms are closed under composition, intersection
and left cancellation. The classes of monomorphisms and of
regular, strong and extremal monomorphisms of $\mathbf{C}$ are
denoted by $\mathit{Mono}$, $\mathit{RegMono}$,
$\mathit{StrongMono}$ and $\mathit{ExtrMono}$, respectively. We
have the inclusions
\[
\mathit{ExtrMono}\subseteq\mathit{StrongMono}
\subseteq\mathit{RegMono}\subseteq\mathit{Mono}.
\]

\begin{theorem}[Theorem 1.3.9 ~\cite{Kurz}]
The following assertions hold in $\mathbf{Set}_F$.

1. (Epi, StrongMono) is a factorization system in
$\mathbf{Set}_F$. Moreover, Epi contains the surjective coalgebra
morphisms, StrongMono contains precisely the injective morphisms,
and (Epi, StrongMono)-factorizations are calculated as
(Epi,Mono)-factorizations in $\mathbf{Set}$.

2. The classes ExtrMono, StrongMono and RegMono coincide.

3. If $F$ preserves weak pullbacks, then Mono coincides with
StrongMono and hence with ExtrMono and RegMono.
\end{theorem}

In virtue of the preceding, unless stated otherwise  we will
assume henceforth that the $\mathbf{Set}$ endofunctor $F$ is
nontrivial, weakly preserves pullbacks and is either bounded or
accessible. Fixing such an $F$, we define the dominion
$\mathcal{M}$ in $\mathbf{Set}_F$ to be the class of monomorphisms
of $\mathbf{Set}_F$  and we set
$\mathbf{Pfn}_F=P(\mathbf{Set}_F,\mathcal{M})$. Next, we extend
the product $\times:\mathbf{Set}_F^2\rightarrow\mathbf{Set}_F$  to
$\boxtimes:\mathbf{Pfn}_F^2\rightarrow\mathbf{Pfn}_F$ by defining
it on objects $X,Y$ by $X\boxtimes Y=X\times Y$ and on partial
morphisms by the formula
\[
\{m,\varphi\}\boxtimes\{m',\varphi'\}=\{m\times
m',\varphi\times\varphi'\}.
\]
This assignment makes sense since $m\times m'\in\mathcal{M}$ and
is independent of the choice of representatives. The reason for
introducing the box product notation is to distinguish it among
the four products in use, namely, the product of categories, the
product in $\mathbf{Set}$, the product in $\mathbf{Set}_F$ and its
extension to a {\it near product}, to be defined, on
$\mathbf{Pfn}_F$. With these definitions,  a coalgebraic version
of Rosolini's Proposition 2.1.1 is immediate ~\cite{Rosolini}.

\begin{proposition}
There is a faithful embedding of $\mathbf{Set}_F$ into
$\mathbf{Pfn}_F$ which sends each $F$-coalgebra  to itself, and
which sends the $F$-coalgebra morphism $\varphi:A\rightarrow B$ to
the equivalence class $\{1_A,\varphi\}$. The product $\times$ on
$\mathbf{Set}_F$ extends to the bifunctor $\boxtimes$ on
$\mathbf{Pfn}_F$.
\end{proposition}

Since the forgetful functor
$U:\mathbf{Set}_F\rightarrow\mathbf{Set}$ creates colimits, it
follows that the category $\mathbf{Pfn}_F$ is cocomplete.

\subsection{$P$- and $P\Sigma$- categories of coalgebras}
The category $\mathbf{Pfn}_F$ with the product $\boxtimes$ has the
structure of a P-category, a notion due to ~\cite{Rosolini} which
generalizes the dominical categories of ~\cite{MR88i:03078} and
the $B$-categories (categories with a binary product) of
~\cite{MR92e:03067}. We recall the relevant definitions from
~\cite{MR92e:03067} with slight changes in notation and with some
additional remarks.

Let $F,G:{\mathbf C}\rightarrow\mathbf{D}$ be functors. An {\it
infranatural transformation} $\phi:F\rightarrow G$ is a family of
morphisms $\phi_X:F X\rightarrow G X$, called components, for each
object $X$ of $\mathbf{C}$. The {\it naturalizer} of $\phi$,
denoted by ${\mathrm nat}\phi$, is the largest subcategory of
$\mathbf{C}$  containing all the objects of $\mathbf{C}$ such that
$\phi$ is a natural transformation $F|\mathrm{nat}\phi\rightarrow
G|\mathrm{nat}\phi$.

Let $\mathbf{C}$ be a category. The {\it diagonal functor}
$\delta_\mathbf{C}:\mathbf{C}\rightarrow\mathbf{C}\times\mathbf{C}$
is given on objects by $A\mapsto(A,A)$ and on morphisms by
$f\mapsto(f,f)$.

A {\it P-category} consists of a category $\mathbf{C}$ together
with a functor
$\boxtimes:\mathbf{C}\times\mathbf{C}\rightarrow\mathbf{C}$,
called a {\it near product}, a natural transformation
$\Delta:1_\mathbf{C}\rightarrow\boxtimes\circ\delta_\mathbf{C}$,
and infranatural transformations $p_i:\boxtimes\rightarrow\pi_i$,
where $\pi_i:\mathbf{C}^2\rightarrow\mathbf{C}$ is the projection
onto the $i$-th factor for $i=0,1$. For morphisms
$f:X\rightarrow Y, g:X\rightarrow Z$ of $\mathbf{C}$ , we set %
$\langle f,\ g \rangle=(f \boxtimes g)\Delta_X :X\rightarrow
Y\boxtimes Z$. These functors and transformations are subject to
the following four conditions:

i) For objects $X$ in $\mathbf{C}$ the following diagrams commute.
\[
\xymatrix{ X\ar[d]_{\Delta_X}\ar[dr]^1 && X\boxtimes
X\ar[r]^-{\Delta_{X\boxtimes X}}\ar[dr]_1 & (X\boxtimes
X)\boxtimes (X\boxtimes X)\ar[d]^-{p_{0 X,X}\boxtimes p_{1 X,X}} \\
X\boxtimes X\ar[r]_-{p_{i X,X}}& X && X\boxtimes X }
\]
where $i=0,1$.

ii) If $\mathbf{P}\subseteq\mathbf{C}$ is the smallest
subcategory closed under $\boxtimes$ containing all components of
$p_i, i=0,1$, then
$\mathbf{C}\times\mathbf{P}\subseteq\mathrm{nat}p_0$ and
$\mathbf{P}\times\mathbf{C}\subseteq\mathrm{nat}p_1$. This
implies that that projections satisfy certain identities; e.g.,
as in the following naturality square diagram, where we have used
the property that the component $p_{0 Y, Z}:Y\boxtimes
Z\rightarrow Y$ is in $\mathbf{P}$, and therefore $(1_X, p_{0
Y,Z})$ is in $\mathrm{nat}p_0$.
\[
\xymatrix{%
(X, Y\boxtimes Z) \ar[d]_{(1_X,p_{0 Y, Z})} & %
X\boxtimes (Y\boxtimes Z) \ar[rr]^-{p_{0 X,Y\boxtimes Z}}
\ar[d]_{1_X\boxtimes p_{0 Y, Z}} && X \ar[d]^{1_X} \\
(X,Y) & X\boxtimes Y \ar[rr]^-{p_{0 X,Y}} && X
}%
\]

iii) There is a natural isomorphism
\[
\xymatrix@R=0.125cm{\mathrm{ass}_\boxtimes:
((-\boxtimes-)\boxtimes-)
\ar[r]& %
(-\boxtimes(-\boxtimes-)) }
\]
of functors $\mathbf{C}^3\rightarrow\mathbf{C}$ whose component
$\mathrm{ass}_{\boxtimes X,Y,Z}$  is given by
\[
\xymatrix@C=0.5cm{{\langle p_{0 X,Y} p_{0 X\boxtimes Y,Z}, %
      \langle p_{1 X,Y} p_{0 X\boxtimes Y,Z}, %
       p_{1 X\boxtimes Y,Z}\rangle\rangle }: (X \boxtimes Y )\boxtimes Z \ar[r]&X
\boxtimes (Y \boxtimes Z). }
\]

iv) Let $\mathrm{tr}_\times$ be the endofunctor on $\mathbf{C}^2$
given by $(X,Y)\mapsto(Y,X)$. There is a natural isomorphism
\[
\xymatrix@R=0.125cm{
\mathrm{tr}_\boxtimes:\boxtimes\ar[r]&\boxtimes\circ\mathrm{tr}_\times
}%
\]
of functors $\mathbf{C}^2\rightarrow\mathbf{C}^2$ whose component
$\mathrm{tr}_{\boxtimes X,Y}$ is given by
\[
\xymatrix@R=0.125cm{
 \langle p_{1 X,Y},\ p_{0 X,Y}\rangle:X\boxtimes Y\ar[r]& Y\boxtimes X
 }.
\]
The natural isomorphisms $\mathrm{ass}_\boxtimes$ and
$\mathrm{tr}_\boxtimes$ must make $\boxtimes$ coherently
associative and commutative; i.e., the natural isomorphism
$\mathrm{tr}_\boxtimes$ must satisfy the condition
$\mathrm{tr}_{\boxtimes Y.X}\circ\mathrm{tr}_{\boxtimes
X,Y}=1_{X\boxtimes Y}$ and a hexagonal coherence condition; the
natural isomorphism $\mathrm{ass}_\boxtimes$ must satisfy a
pentagonal coherence condition ~\cite{MR1712872,MR30:1160}.

\begin{proposition}
The category $\mathbf{Pfn}_F$ is a $P$-category.
\end{proposition}
\begin{proof}
This follows from Rosolini's Theorem 2.1.9 ~\cite{Rosolini}
applied to $\mathbf{Set}_F$.
\end{proof}

\subsection{Prodominical categories}
A {\it system of zero morphisms} is a collection of morphisms
$0_{X,Y}:X\rightarrow Y$  for each pair of objects $X$ and $Y$ of
${\mathbf C}$ such that for objects $W, Z$ and morphisms
$f:W\rightarrow X$ and $g:Y\rightarrow Z$ of $\mathbf{C}$, one has
$g 0_{X,Y}f = 0_{W,Z}$. A system of zero morphisms is unique if it
exists. For any coalgebra $A$, the empty map $0_{\emptyset,A}:
\emptyset\rightarrow A$ is vacuously a morphism in
$\mathbf{Set}_F$. The equivalence classes
$\{0_{\emptyset,A},0_{\emptyset,B}\}$ for coalgebras $A$ and $B$
form a system of zeros of $\mathbf{Pfn}_F$.

A {\it prodominical} category $\mathbf{C}$ is a $P$-category  that
is {\it pointed}; i.e., $\mathbf{C}$ contains a system of zero
morphisms and, for any $\phi:A\rightarrow B$,
$\phi\boxtimes0_{C,D}=0_{A\boxtimes C,B\boxtimes D}$. A morphism
$f:X\rightarrow Y$ of a pointed category is {\it weakly total} if
for all $\phi:W\rightarrow X$, $f\phi=0_{W,Y}$ implies that
$\phi=0_{W,X}$.

We take the following proposition from Heller ~\cite{MR92e:03067}
as the definition of a dominical category. A prodominical category
is {\it dominical} if every weakly total morphism is total and if
$\phi\boxtimes\psi=0$ implies $\phi=0$ or $\psi=0$. This property
fails in $\mathbf{Pfn}_F$ in general, for two reasons.

The category $\mathbf{Pfn}_F$ can fail to be dominical due to
certain pathologies of the product in $\mathbf{Set}_F$; Gumm and
Schr{\"o}der ~\cite{GS2000} give examples of three finite
coalgebras $A,B,C$ for the finite powerset functor such that
$A\times A\cong A$, $A\times B =\emptyset$ and such that $C\times
C$ is infinite.

\begin{proposition} There exists a nontrivial, bounded
$\mathbf{Set}$ endofunctor $F$ that weakly preserves pullbacks
such that $\mathbf{Pfn}_F$ is prodominical but not dominical.
\end{proposition}
\begin{proof}
The finite powerset functor $P_\omega$ is nontrivial, preserves
weak pullbacks and is bounded
~\cite{GS2000,MR1787576,MR99h:18012}, hence
$\mathbf{Pfn}_{P_\omega}$ is prodominical. Moreover there exist
two nonempty finite coalgebras $A$ and $B$ for the finite powerset
functor such that $A\times B=\emptyset$ ~\cite{GS2000}. Let
$i_A:A\rightarrow A\coprod B$ be the coproduct injection in
$\mathbf{Set}_{P_\omega}$, and let $1_B$ be the identity on $B$.
By previous remarks,
\[
i_A\times1_B:A\times B\rightarrow (A\coprod B)\times B =
0_{\emptyset,B\times B},
\]
and therefore there exist nonzero
$\phi,\psi\in\mathbf{Pfn}_{P_\omega}$ with $\phi\boxtimes\psi=0$.
\end{proof}

Also, the category $\mathbf{Pfn}_F$ can fail to be dominical if it
contains a non-total weakly total morphism. An example of a
prodominical category that is not dominical because it has such a
morphism is the syntactic $P$-category $T$, obtained from Peano
arithmetic {\sl PA}; this example is essentially due to Franco
Montagna ~\cite{Rosolini}. The category $T$ has exactly one
object. The morphisms of $T$ are the equivalence classes under
provability of of the provably functional $\Sigma_1$-formula
$F(x,y)$ of {\sl PA}; i.e., those $\Sigma_1$ $F(x,y)$ for which
\[
\vdash_{\sl PA} \forall x, y, z (F(x,y)\wedge F(x,z)\Rightarrow
y=z).
\]
The identity morphism $1_T$ is defined to be the provable
equivalence class of the formula $0=0$. The composition of
morphisms $F(x,y)$ and $G(y,z)$ is defined by $(G\circ
F)(x,z)\Leftrightarrow \exists y F(x,y)\wedge G(y,z)$. The domain
of a morphism $F(x,y)$ is the morphism $D(x,y)$ defined
by $\exists z F(x,z)\wedge y=z$. A morphism is total if %
$\vdash_{\sl PA}\forall x D(x,x)$. A morphism $F(x,y)$ is a zero
morphism (undefined) provided $\vdash_{PA}\forall x,y \lnot
F(x,y)$. A morphism $F(x,y)$ is weakly total if for all morphisms
$G(w,x)$,
\[
\vdash_{PA}  \forall w,y \lnot(F\circ G)(w,y)  \text{\ implies\ }
\vdash_{PA}  \forall w,x \lnot(F\circ G)(w,x).
\]
See Rosolini ~\cite{Rosolini} for the argument that Peano
arithmentic contains a provably functional $\Sigma_1$ formula
$F(x,y)$ that is weakly total but not total. Further examples of
prodominical categories that are not dominical are given in
Montagna ~\cite{MR90i:03046}.

However, it is even easier to produce a category of coalgebras and
partial morphisms which possesses a non-total weakly total
morphism.  In $\mathbf{Pfn}_F$, weak totality can be interpreted
by means of the coalgebra topology. The topological criterion
given below for a weakly total morphism to be total  can be used
to give an elementary proof that $\mathbf{Pfn}_F$ is not dominical
in general. The following statements are immediate.

\begin{proposition} \label{PropTopology}
Let $\varphi:X\rightarrow Y = \{m:W\rightarrowtail X,\psi:W\rightarrow Y\}$ %
be a morphism in $\mathbf{Pfn}_F$, with $m$ mono.

i) The morphism $\varphi$ is weakly total if and only if the image
$m[W]$ of $m$ is dense in the coalgebra topology on $X$.

ii) If the subcoalgebra $U$ of $X$ is dense in $X$, then the
partial morphism $\{i, 1_U\}$ defined by
\[
\xymatrix{{\ U}\ar@{^{(}->}[d]_i\ar[r]^{1_U} & U\\
X \ar@{-->}[ur]_-{\{i,1_U\}} &}
\]
is weakly total.

iii) The coalgebra $X$ is irreducible in the coalgebra topology if
and only if every morphism $\varphi:X\rightarrow Y$ other than
$0_{X,Y}$ is weakly total.

iv) Every weakly total morphism $\varphi:X\rightarrow Y$ is total
if and only if $X$ is the only nonempty open dense subset of $X$.

v) If $X$ is Hausdorff, then every weakly total morphism
$\varphi:X\rightarrow Y$ is total.
\end{proposition}

To obtain a non-total weakly total morphism in $\mathbf{Pfn}_F$,
it suffices to find a coalgebra with a proper nonempty
subcoalgebra that is dense in the coalgebra topology.

\begin{proposition} Let $1_\mathbf{Set}$ be the identity functor on $\mathbf{Set}$.
Then $\mathbf{Pfn}_{1_\mathbf{Set}}$ is prodominical but not
dominical.
\end{proposition}
\begin{proof}  A coalgebra for the identity functor
is given by a self map of a set. Let $X=\{x, y\}$ be a set with
$x\ne y$, and define $\alpha:X\rightarrow X$ by $\alpha(x) = y$
and $\alpha(y)=y$. Then $U=\{y\}$ is a proper subcoalgebra of $X$
that is open and dense in the coalgebra topology on $X$. It
follows from Proposition \ref{PropTopology}, statement iv)  that
the partial map $X\rightarrow U$ of statement ii) is weakly total
but not total.
\end{proof}

We record some non-pathological properties of the product in
$\mathbf{Set}_F$ that will be needed in the sequel. Fortunately,
despite the example of two finite nonempty coalgebras $A,B$ for
which $A\times B=\emptyset$, we have the following.

\begin{proposition}\label{PropPower}
If $\mathbf{Set}_F$ has finite products, and if $A$ is a nonempty
$F$-coalgebra, then for each $n\ge 1$, $A^n$ is nonempty.
\end{proposition}
\begin{proof}
Let $U:\mathbf{Set}_F\rightarrow\mathbf{Set}$ be the forgetful
functor. If $A$ is nonempty, so is the diagonal
$\Delta_A=\{(a,a):a\in U(A)\}$, which is a bisimulation on $A$
~\cite{MR1791953}.  Lemma 8.1 of ~\cite{GS2000} states that if
$A_0\times A_1$ is a product in $\mathbf{Set}_F$ with projections
$\pi_i, i=0,1$, then $\left(U\pi_0\times U\pi_1\right)U(A_0\times
A_1)$ is the largest bisimulation between $A_0$ and $A_1$; taking
$A_0=A_1=A$, we have that $\emptyset\ne \Delta_A\subseteq
\left(U\pi_0\times U\pi_1\right)U(A\times A)$ so that $A\times A$
cannot be empty. It follows that $A^n\ne\emptyset$ for any $n\ge
1$.
\end{proof}

The product in $\mathbf{Set}_F$ distributes over the
coproduct.\footnote{Necessary and sufficient conditions for
distributivity in categories of coalgebras were announced in joint
work of H. Peter Gumm and Tobias Schr\"oder of the University of
Marberg, Germany, and Jesse Hughes of Carnegie Mellon University
~\cite{GHSdist}.} This follows essentially from theorem 6.4 of
Worrel ~\cite{MR99k:18007}, which states that if $F$ is a bounded
endofunctor which preserves weak pullbacks, then $\mathbf{Set}_F$
is a full reflective subcategory of a Grothendieck topos.
Moreover, assuming only that $F$ is nontrivial and preserves weak
pullbacks, Worrel's proof implies that the following conditions of
Giraud's theorem hold. In $\mathbf{Set}_F$, coproducts
 are disjoint and stable under pullback, every epimorphism is the
 coequalizer of its kernel pair, and epimorphisms are stable under
 pullback.

\begin{proposition}\label{PropDist}
Let $F$ be a nontrivial $\mathbf{Set}$ endofunctor that weakly
preserves pullbacks and such that $\mathbf{Set}_F$ has finite
limits. Let $X$ be an $F$-coalgebra, let $I$ be a set, and let
$\{Y_\alpha\}_{\alpha\in I}$ be a family of $F$-coalgebras. There
is a canonical isomorphism
\[
\xymatrix{ X\times\coprod_I Y_\alpha \ar[r]&\coprod_I \left(X\times Y_\alpha\right) } %
\]
\end{proposition}
\begin{proof}
By a result of Gumm and Schr\"oder ~\cite{MR1787576}, the preimage
of a subcoalgebra is a subcoalgebra if and only if $F$ weakly
preserves pullbacks along monomorphisms; hence pullbacks of
subcoalgebras are subcoalgebras in $\mathbf{Set}_F$. Paraphrasing
Worrel ~\cite{MR99k:18007}, since the forgetful functor
$U:\mathbf{Set}_F\rightarrow \mathbf{Set}$ creates colimits and
since $U$ weakly preserves pullbacks if $F$ does, coproducts are
disjoint and stable under pullback. Fix coalgebras $A$ and $B$; it
follows by standard arguments ~\cite{SMIM} that the pullback
operation $-\times_A B$ distributes over coproducts; taking $A$ to
be the terminal object, it follows that products distribute over
coproducts.
\end{proof}

In a $P$-category $\mathbf{C}$, the {\it domain}
$\mathrm{dom}\phi$ of a morphism $\phi:X\rightarrow Y$
is the composite %
$p_{0 X,Y}\circ\langle 1_X\ \phi\rangle:X\rightarrow X$. For an
object $X$ of $\mathbf{C}$, let $\mathrm{dom}(X)$ denote the set
of domains $\mathrm{dom}\phi$ for morphisms $\phi:X\rightarrow Y$.
The set $\mathrm{dom}(X)$ is a meet semilattice with meet defined
by composition ~\cite{MR88i:03078,MR92e:03067,Rosolini}. Domains
in $\mathbf{Pfn}_F$ correspond precisely to subcoalgebras.

\begin{proposition} \label{PropMeetSemi}
Let $X$ be a coalgebra in $\mathbf{Pfn}_F$. There is a meet
semilattice isomorphism from $\mathrm{dom}(X)$ to the lattice
$\mathcal{L_X}$ of sub-coalgebras of $X$ (considered as a meet
semilattice).
\end{proposition}
\begin{proof}
The map $\mathcal{L_X}\rightarrow\mathrm{dom}(X)$ is defined by
sending the sub-coalgebra $U$ of $X$ to the domain $p_{0
X,Y}\circ\langle 1_X\ \phi\rangle$, where $\phi=\{U\hookrightarrow
X,U\hookrightarrow X\}$. The inverse map is defined as follows.
Given a domain $\phi\in\mathrm{dom}(X)$, we may write %
$\phi= p_{0 X,Y}\circ\langle 1_X\ \{m,\varphi\}\rangle$, where
$m:V\rightarrow X$ is a monomorphism of the coalgebras, and where
$\varphi:V\rightarrow Y$ is a morphism in $\mathbf{Set}_F$. By
definition of composition in $\mathbf{Pfn}_F$ one obtains the
following diagram in which the squares are pullbacks. We adopt the
convention of not showing the structure maps and induced morphisms
of the coalgebras occurring in commutative diagrams of coalgebras.
\[
\newdir{ >}{{}*!/-10pt/@{>}}
\xymatrix{
V\ar[r]\ar[d]^1 %
& X\times V\ar[r]^-{1_X\times \varphi}\ar[d]^1
& X \times Y\ar[d]_1\ar[r]^-{p_{0 X,Y}} %
& X %
\\ %
V \ar[r]\ar@{ >->}[d] %
& X\times V\ar@{ >->}[d]\ar[r]^-{1_X\times \varphi} %
& X\times Y\ar@{-->}[ur]_{p_{0 X,Y}} %
\\ %
X \ar[r]^-{\Delta_X}\ar[d]^1 %
& X\times X\ar@{-->}[ur]_{1_X\boxtimes\varphi} %
\\ %
X \ar@{-->}[ur]_{\Delta_X} }
\]
We use this diagram to define the coalgebra $U\in\mathcal{L_X}$
corresponding to the domain $\phi$ as the image of $V$ under the
left vertical column of the diagram. Showing that the meets are
preserved involves certain large diagrams such as those occurring
in ~\cite{MR88i:03078}.
\end{proof}

Let $\mathbf{C}$ be  a $P$-category. A morphism $\phi:X\rightarrow
Y$ $\mathbf{C}$ is {\it total} if $\mathrm{dom}\phi=1_X$. The
collection of total morphisms of $\mathbf{C}$ generate its
subcategory $\mathbf{C}_T$ of total morphisms. The near-product
and infranatural transformations of a $P$-category $\mathbf{C}$
become a product and natural transformations, respectively, on its
subcategory $\mathbf{C}_T$, which has the structure of a
$B$-category; for convenience we include the definition from
Heller ~\cite{MR92e:03067}. A {\it $B$-category} is a category
$\mathbf{C}$ with a bifunctor
$\times:\mathbf{C}^2\rightarrow\mathbf{C}$ and natural
transformations
$p_0,p_1,\Delta,\mathrm{ass}_\times,\mathrm{tr}_\times$ satisfying
the conditions i), iii) and iv) for $P$-categories above, with
$\times$ replacing $\boxtimes$. It should be emphasized that the
projections of a $B$-category are required to be natural and not
merely infranatural transformations; we write $\times$ for the
product of a $B$-category and $\boxtimes$ for the near product of
a $P$-category. Dually, one may speak of a category with a binary
coproduct, together with natural
transformations %
$i_0,i_1,\nabla,\mathrm{ass}_{\coprod},\mathrm{tr}_{\coprod}$
satisfying the duals of the conditions 1), iii) and iv) in which
the {\it injections} $i_0,i_1$ replace the {\it projections}
$p_0,p_1$, the {\it codiagonal} $\nabla$ replaces the {\it
diagonal} $\Delta$, where $\mathrm{ass}_{\coprod}$ and
$\mathrm{tr}_{\coprod}$ replace $\mathrm{ass}_\times$ and
$\mathrm{tr}_\times$, respectively, and where $\coprod$ replaces
$\boxtimes$. A $B$-category with a binary coproduct that has a
natural isomorphism called $\mathit{dist}$ inverse to the natural
transformation
\[
\left(X\times Y \right)\coprod \left(X\times Z \right)\rightarrow
X \times\left(Y\coprod Z\right)
\]
is called a {\it $B^+$-category}.  A $B^+$-category with a
countable coproduct $\coprod_\mathbb{N}$ such that $\times$
distributes over $\coprod_\mathbb{N}$ is called a %
{\it $B\Sigma$-category}.

 By analogy with $B^+$- and $B\Sigma$-categories,
one may define $P^+$ and $P\Sigma$-categories, in which the
coproduct injections are required to be natural (and not merely
infranatural) transformations. Under the embedding
$\mathbf{Set}_F\rightarrow\mathbf{Pfn}_F$, the canonical
isomorphism of proposition \ref{PropDist} defines the natural
isomorphism
\[
\xymatrix{ \mathit{dist}_{X,(Y_n)}:X\boxtimes\coprod_\mathbb{N}
Y_n\ar[r]&\coprod_\mathbb{N}X\boxtimes Y_n }
\]
 of functors
$\mathbf{Pfn}_F\times\mathbf{Pfn}_F^{\mathbb{N}}\rightarrow
\mathbf{Pfn}_F$.

\begin{proposition}
The category $\mathbf{Pfn}_F$ is a $P\Sigma$-category.
\end{proposition}
\begin{proof}
Let $f:X\dashrightarrow W$ in $\mathbf{Pfn}_F$ be given by
$f=\{\mu,\phi\}$, where $\mu:U\rightarrowtail X$ and
$\phi:U\rightarrow W$ are in $\mathbf{Set}_F$, with $\mu$ mono,
and for $n\in\mathbb{N}$ let $g_n:Y_n\dashrightarrow Z_n$ in
$\mathbf{Pfn}_F$ be given by  $g_n=\{\nu_n,\psi_n\}$, where
$\nu_n:V_n\rightarrowtail Y_n$ and $\psi_n:V_n\rightarrow Z_n$ are
in $\mathbf{Set}_F$, with $\nu_n$ mono. Consider the following
diagram.
\[
\newdir{ >}{{}*!/-10pt/@{>}}
\xymatrix{ U\times\coprod_\mathbb{N} V_n %
   \ar[rr]^-{\mathit{dist}_{U,(V_n)} } %
   \ar@{ >->}[dr]^-{\mu\times\coprod_\mathbb{N}\nu_n}
   \ar[dd]^{\phi\times\coprod_\mathbb{N}\psi_n}
                    && \coprod_\mathbb{N}(U\times V_n)  %
                       \ar@{ >->}[dr]^-{\coprod_\mathbb{N}(\mu\times\nu_n)} %
                       \ar'[d]^(.7){\coprod_\mathbb{N}(\phi\times\psi_n)}[dd]  %
\\
&  X\times\coprod_\mathbb{N} Y_n %
   \ar[rr]^(.35){\mathit{dist}_{X,(Y_n)}} %
   \ar@{-->}[dl]^{f\boxtimes\coprod_\mathbb{N}g_n} %
                          & { } & \coprod_\mathbb{N} (X\times Y_n) %
                          \ar@{-->}[dl]^{\coprod_\mathbb{N}(f\boxtimes g_n)} %
\\
 W\times\coprod_\mathbb{N} Z_n \ar[rr]_-{\mathit{dist}_{W,(Z_n)}}
                     && \coprod_\mathbb{N}(W\times Z_n)  %
}
\]
We claim that the bottom parallelogram commutes in
$\mathbf{Pfn}_F$. Observe that the top parallelogram is a pullback
in $\mathbf{Set}_F$. It follows that the composite
\[
\coprod_\mathbb{N}\left(f\boxtimes g_n\right)\circ\mathit{dist}_{X,(Y_n)} %
\]
in $\mathbf{Pfn}_F$ is represented by the pair
\begin{align}\label{EqnPair}
\left\{\mu\times\coprod_\mathbb{N}\nu_n,
\coprod_\mathbb{N}\left(\phi\times\psi_n\right)\circ\mathit{dist}_{U,(V_n)}\right\}.
\end{align}
On the other hand, the composite
\[
\mathit{dist}_{W,(Z_n)}\circ\left(f\boxtimes\coprod_\mathbb{N}
g_n\right)
\]
in $\mathbf{Pfn}_F$ is represented by the pair
\[
\left\{\mu\times\coprod_\mathbb{N}\nu_n, %
\mathit{dist}_{W,(Z_n)}\circ\left(\phi\times\coprod_\mathbb{N}\psi_n\right)\right\}, %
\]
which equals (\ref{EqnPair}) since the (back) rectangle commutes,
as {\it dist} is natural in $\mathbf{Set}_F$.
\end{proof}

In a $P$-category, if $\phi:X\rightarrow Y$ is a morphism and if
$\varepsilon\in\mathrm{Dom}Y$, we write $\phi\prec\varepsilon$ if
and only if $\varepsilon\phi=\phi$; we say that $\varepsilon$ {\it
receives} $\phi$. If $\varepsilon$ receives $\phi$, and in
addition, $\varepsilon$ satisfies for all appropriate
$\psi,\psi'$, $\psi\phi=\psi'\phi$ implies
$\psi\varepsilon=\psi'\varepsilon,$ then $\varepsilon$ is the
least domain in $\mathrm{Dom}Y$ receiving $\phi$, since if
$\delta\in\mathrm{dom}Y$ satisfies $\delta\phi=\phi$, then
$\delta\phi=\varepsilon\phi$, which implies that
$\delta\varepsilon=\varepsilon\varepsilon=\varepsilon$ and
therefore $\varepsilon\prec\delta$. In this case we say that
$\phi$ has range $\varepsilon$ and we write
$\mathrm{ran}\phi=\varepsilon$.

If each morphism of the $P$-category $\mathbf{C}$ has a range, we
say that $\mathbf{C}$ {\it has ranges}. Since $\mathbf{Set}_F$ has
(Epi,StrongMono)-factorizations, the image of a coalgebra under a
morphism is a subcoalgebra of the codomain; conversely every
subcoalgebra is an image. If $\phi=\{m,\psi\}$ is a morphism of
$\mathbf{Pfn}_F$, where $m:U\rightarrow X$ and $\psi:U\rightarrow
Y$ are morphisms of $\mathbf{Set}_F$ with $m$ mono, we define the
image $\mathrm{im}\phi$ of $\phi$ by $\mathrm{im}\phi=\psi[U]$.
Under the semilattice isomorphism of Proposition
\ref{PropMeetSemi}, images of morphisms in $\mathbf{Pfn}_F$, which
are coalgebras, correspond with ranges, which are morphisms.

\begin{proposition}
The category $\mathbf{Pfn}_F$ has ranges.
\end{proposition}

If every morphism of  the $P$-category $\mathbf{C}$ has a range,
and if for morphisms $\phi,\psi$ of $\mathbf{C}$,
$\mathrm{ran}(\phi\boxtimes\psi)=\mathrm{ran}\phi\boxtimes\mathrm{ran}\psi$,
then one says ~\cite{MR92e:03067} that $\mathbf{C}$ is an {\it
$rP$-category}; such a category has a calculus of ranges
~\cite{MR88i:03078}. In a $B^+$-category, if $f$ and $g$ are
morphisms with the same codomain $Y$, then we define $[f,\
g]=\nabla_Y(f\coprod g)$. In an $rP^+$- ($rP\Sigma$-) category,
the meet semilattice $\mathrm{dom}(X)$ becomes a distributive
lattice if one defines the join by
$\varepsilon\cup\delta=\mathrm{ran}[\varepsilon,\ \delta]$ for
$\varepsilon,\delta\in\mathrm{dom}(X)$
~\cite{MR88i:03078,MR92e:03067}. To show that $\mathbf{Pfn}_F$ is
an $rP\Sigma$-category, we first show that the image of a product
of morphisms in $\mathbf{Set}_F$ is the product of the images; the
proof is mostly folklore. The lack of a simple description of the
product in $\mathbf{Set}_F$ and the possibility of pathologies
seems to necessitate the use of categorical technique in the
proof.

\begin{proposition}
For morphisms $\phi,\psi$ of $\mathbf{Set}_F$,
\[
\mathrm{im}(\phi\times\psi)=\mathrm{im}\phi\times\mathrm{im}\psi.
\]
\end{proposition}
\begin{proof}
We first show that if $f:A\rightarrow B$ is an epimorphism in
$\mathbf{Set}_F$, then so is $f\times 1:A\times C\rightarrow
B\times C$. The following diagram is a pullback in
$\mathbf{Set}_F$.
\[
\xymatrix{ A\times C \ar[r]^-{f\times 1}\ar[d]_{p_A} & B\times
C\ar[d]^-{p_B} \\
A \ar[r]^-f & B}
\]

Apply the forgetful functor $U:\mathbf{Set}_F\rightarrow
\mathbf{Set}$ to this and take the pullback in $\mathbf{Set}$.
Proposition 6.3 of ~\cite{MR99k:18007} states that if $F$
preserves weak pullbacks, then the forgetful functor $U$ preserves
weak pullbacks; hence the square in the following diagram is a
weak pullback, and therefore there exists an induced map as
indicated.
\[
\xymatrix{U A \times_{U B} U(B\times C) \ar@/^2pc/[drrr]\ar@{.>}[dr]\ar@/_1pc/[ddr] & &  &\\
& U(A\times C) \ar[rr]^-{U(f\times 1)}\ar[d]_{U p_A} && U(B\times
C)\ar[d]^-{U p_B} \\
& U A \ar[rr]^-{U f} && U B }
\]
Since $U$ preserves epimorphisms, $U f$ is surjective. Since it is
possible that $U(B\times C)=\emptyset$, in that trivial case
commutativity of the square forces $U(A\times C)=\emptyset$, so
that $U(f\times1)=0_{\emptyset,\emptyset}=1_\emptyset$. Otherwise,
a diagram chase shows that $U(f\times 1)$ is surjective and, since
$U$ reflects epimorphisms, $f\times1$ is an epimorphism in
$\mathbf{Set}_F$. Therefore, if $f,g$ are epimorphisms in
$\mathbf{Set}_F$, so is $f\times g= (f\times 1)\circ(1\times g)$.

Using the first pullback diagram alone, one has that if $f$ is a
monomorphism, then so is $f\times1$, and therefore a product of
monomorphisms in $\mathbf{Set}_F$ is a monomorphism. It follows
that the epi-mono factorization of a product $f\times g$ is the
product of the epi-mono factorizations of $f$ and of $g$.

A theorem of Gumm and Schr\"oder ~\cite{GS2000} states that if
$A_i$ is a subcoalgebra of $B_i$ for $i=0,1$, and if the product
$B_0\times B_1$ exists, then so does $A_0\times A_1$, which is a
subcoalgebra of $B_0\times B_1$. If $\phi:X\rightarrow Y$ and
$\psi:W\rightarrow Z$ are coalgebra morphisms, then by previous
remarks, one has two canonically isomorphic epi-mono
factorizations of $\phi\times\psi$ as indicated in the following
commutative diagram.
\[
\xymatrix{ X\times W\ar@{>>}[r]\ar@{>>}[d]& \mathrm{im}(\phi\times\psi)\ar@{^{(}->}[d] & \\
 \mathrm{im}\phi\times\mathrm{im}\psi\ar@{^{(}->}[r]\ar@{.>}[ur] & Y\times Z}
\]
Since the morphisms into $Y\times Z$ are inclusions, the induced
map is the identity.
\end{proof}

\begin{proposition}
The category $\mathbf{Pfn}_F$ is an $rP\Sigma$-category. Moreover,
for any $F$-coalgebra $X$, the meet semilattice isomorphism of
Proposition \ref{PropMeetSemi} is a lattice isomorphism from
$\mathrm{dom}(X)$ to the lattice $\mathcal{L_X}$ of sub-coalgebras
of $X$.
\end{proposition}

The following definitions are taken from Heller
~\cite{MR92e:03067}. In a $P$-category, a {\it section} of a
morphism $\phi:X\rightarrow Y$  is a morphism
$\sigma:Y\rightarrow X$ such that $\phi\sigma=\mathrm{dom}\sigma$
and $\phi\sigma\phi=\phi$. A $P$-category satisfies the {\it
axiom of choice} if every morphism has a section. A morphism
$\phi$ is a {\it partial monomorphism} if
$\phi\theta=\phi\theta'$ implies
$(\mathrm{dom}\phi)\theta=(\mathrm{dom}\phi)\theta'$. A
$P$-category satisfies the {\it weak axiom of choice} if each
partial monomorphism has a section.

\begin{proposition}
The category $\mathbf{Pfn}_F$ satisfies the weak axiom of choice.
\end{proposition}
\begin{proof}
Let $\phi:X\rightarrow Y$ be a partial monomorphism in
$\mathbf{Pfn}_F$. Then  $\phi=\{\mu,\psi\}$, where $\mu, \psi
:U\rightarrow X$ are monomorphisms in $\mathbf{Set}_F$, and where
$U$ is a subcoalgebra of $X$. The category $\mathbf{Set}_F$ has
(Epi-StrongMono) factorizations, hence one has the following
diagram in which $\psi$ has the epi-mono factorization $\psi=m
e$. \[
\newdir{ >}{{}*!/-10pt/@{>}}
\xymatrix{ %
Q \ar@<1ex>@{ >->}[d]^f\ar@{ >->}[r]^m & Y \ar[d]^1 \\ %
U \ar@<1ex>@{>>}[u]^e\ar@{ >->}[d]^\mu\ar@{ >->}[r]^-{\psi} & Y \\ %
X\ar@{-->}_-{\phi=\{\mu, \psi\}}[ur] }
\]
A morphism in $\mathbf{Set}_F$ is epi if and only if it is
surjective in $\mathbf{Set}$. The morphism $e$ is a monomorphism
since $m$ and $\psi$ are. Therefore $e$ is invertible in
$\mathbf{Set}$ with inverse $f$ which, since the image of $m$ is
contained in the image of $\psi$, is a morphism in
$\mathbf{Set}_F$ by a diagram lemma (e.g., Gumm ~\cite{GS2000}).

The partial morphism $\sigma$ defined by $\sigma=\{m,\mu f\}$
gives the required section. The required properties
$\phi\sigma=\mathrm{dom}\sigma$ and $\phi\sigma\phi=\phi$ can be
verified in the following commutative diagram, in which the
squares are pullbacks in $\mathbf{Set}_F$.
\[
\newdir{ >}{{}*!/-10pt/@{>}}
\xymatrix{
U\ar@{>>}[r]^e\ar[d]^1 %
& Q\ar@{ >->}[r]^-{f}\ar[d]^1
& U \ar@{ >->}[d]^\mu\ar[r]^-{\psi} %
& Y %
\\ %
U \ar@{>>}[r]^e\ar[d]^1 %
& Q \ar@{ >->}[d]^m\ar@{ >->}[r]^-{\mu f} %
& X\ar@{-->}[ur]_{\phi} %
\\ %
U \ar[r]^-{\psi}\ar@{ >->}[d]^\mu %
& Y\ar@{-->}[ur]_{\sigma} %
\\ %
X \ar@{-->}[ur]_{\phi} }
\]
\end{proof}

In consequence, a generalization of what Gumm ~\cite{glutacs}
refers to as the first diagram lemma of $\mathbf{Set}_F$ holds in
$\mathbf{Pfn}_F$. In Heller's nomenclature, one says that
$\mathbf{Set}_F$ and $\mathbf{Pfn}_F$ are {\it factorial}
~\cite{MR92e:03067}.

\begin{lemma}
Let $X,Y$ and $Z$ be coalgebras, let $\phi:X\rightarrow Y$ and
$\psi:X\rightarrow Z$ be partial morphisms in $\mathbf{Pfn}_F$
with $\mathrm{dom}\phi\subseteq\mathrm{dom}\psi$, and suppose that
for each $x,y\in\mathrm{dom}\phi$, $\phi x = \phi y$ implies that
$\psi x = \psi y$ (in Heller's terminology one says that $\phi$
{\it divides} $\psi$). Then there is a unique partial map, denoted
by $\psi/\phi$, such that $\mathrm{dom}\psi/\phi=\mathrm{ran}\phi$
and such that the following diagram commutes.
\[
\xymatrix{ X \ar[r]^-{\phi} %
\ar[dr]_-{\psi} & Y
\ar@{-->}[d]^-{\psi/\phi} %
\\
& Z }
\]
\end{lemma}

\subsection{Free semigroups, $B^\#$- and iteration categories}

If $\mathbf{C}$ is a $B\Sigma$- or $P\Sigma$-category, one may
construct the {\it free semigroup}
$X^\#=\coprod_\mathbb{N}X^{n+1}$ of an object $X$ of $\mathbf{C}$.
In view of the possibility of pathological products of coalgebras,
we remark that Proposition \ref{PropPower} implies that the free
semigroup of a nonempty coalgebra is nonempty. The associative
multiplication $m:X^\#\times X^\#\rightarrow X^\#$ of the free
semigroup is the composite
\[
\xymatrix{
\left(\coprod_\mathbb{N}X^{n+1}\right)\times\left(\coprod_\mathbb{N}X^{n+1}\right) %
\ar[r] & %
\coprod_{\mathbb{N}\times\mathbb{N}}X^{n+1}\times X^{m+1} %
\ar[r] & \coprod_\mathbb{N}X^{k+1}}
\]
where the left map is obtained from two applications of the
natural transformation $\mathit{dist}$, and where the right map
is, to within a natural isomorphism, the coproduct of the
codiagonal maps
\[
\xymatrix{ \nabla_k:\coprod_{n+m+1=k} X^{n+1}\times X^{m+1} \ar[r] & %
X^{k+1}}
\]
for $k\ge 1$ (compare with the construction of the free monoid in
~\cite{MR1712872}, page 172). In particular, the categories
$\mathbf{Set}_F$ and $\mathbf{Pfn}_F$ possess these canonical
coalgebra morphisms.

More general than $B\Sigma$-categories are the $B^\#$-categories,
whose definition will be needed for the statement of Heller's
existence theorem. Though our results could be formulated in
certain $P\Sigma$-categories without mention of $B^\#$-categories
and iteration categories (to be defined), technicalities involving
products of coalgebras arise that compel us to consider the more
general notion even when the free semigroup is available. A
$B^\#$-category has an associated monad, denoted by $\#$, called
the {\it formally-free semigroup} functor which, although lacking
the universal property of the free semigroup, has enough of the
structure of the free semigroup to define by categorical means
operations corresponding to the manipulation of finite sequences.
Such operations are used to define the iteration of a morphism in
the sequel.

We remark on certain technicalities of products of coalgebras
alluded to above. In the case of coalgebras for an endofunctor $F$
which preserves products, the product of coalgebras may be
constructed as it is constructed in $\mathbf{Set}$, and one may
identify an element of a finite product of $F$-coalgebras with a
finite sequence of elements of coalgebras. However, in most cases
of interest in computer science, there is no simple description of
products of coalgebras; products are not constructed from products
of base sets.  For the proof of the iteration lemma below, one
must rely on the operations provided by the (formally) free
semigroup in the general case.

The definition of a $B^\#$-category, taken verbatim from Heller,
follows. For this purpose some additional notation is needed. The
formally-free semigroup functor $\#$ comes with an associative
multiplication $m$; for morphisms $f,g$ with common codomain
$X^\#$, one defines $f\cdot g= m(f\times g)$. We define
$j_n:X^n\rightarrow X^\#$ by $j_1=j$ and $j_{n+1}=j\cdot j_n$ for
$j>1$. When $X^\#$ is the free semigroup, the $j_n$ are (to within
canonical isomorphism) the coproduct injections.

A $B^\#$-category is a $B^+$-category equipped with seven natural
transformations
\[
\xymatrix@R=0.125cm{
 m:X^\#\times X^\# \ar[r] & X^\# \\
j:X\ar[r] & X^\# \\
e:X^{\#\#}\ar[r] & X^\# \\
l: X^\#\ar[r] & X \coprod \left(X\times X^\#\right) \\
r: X^\#\ar[r] & \left(X^\#\times X\right) \coprod X \\
{\it par}:X\times Y^\#\ar[r]&\left(X\times Y\right)^\#}
\]
\[
\xymatrix@R=0.125cm{
{\it wd}: \left(X\coprod Y\right)^\#\ar[r] & %
{X^\# \coprod Y^\#} \coprod \left(X^\#\times Y^\#\right)^\# \coprod %
\left(Y^\#\times X^\#\right)^\# \\
&\coprod \left(Y^\#\times\left(X^\#\times Y^\#\right)^\# \right)\coprod %
\left(X^\#\times\left(Y^\#\times X^\#\right)^\#\right) }
\]
subject to the following conditions.

i) $m$ is associative, hence $X^\#$ is a semigroup and, if $f$ is
a morphism, then $f^\#$ is a semigroup homomorphism.

ii) $\#$ is a monad with unit $j$ and multiplication $e$, and each
component of $e$ is a semigroup homomorphism.

iii) $l$ and $r$ are respectively the inverses of
\[
\xymatrix@R=0.0625cm{ [j, j\cdot X^\#]:X\coprod \left(X\times X^\#\right)\ar[r] &X^\#, \\ %
[X^\#\cdot j, j]:\left(X^\#\times X\right)\coprod X\ar[r] & X^\#.
}
\]

iv) The following diagram commutes.
\[
\xymatrix{ X\times
Y^\#\ar[r]^-{\mathit{par}}\ar[dr]_-{p_{1,X,Y^\#}} &
\left(X\times Y\right)^\# \ar[d]^{p^{\#}_{1,X,Y}} \\
& Y^\# }
\]

v) $\mathit{wd}$ is inverse to the morphism
\[
[i^\#_0, i^\#_1, e(i^\#_0\cdot i^\#_1)^\#, e(i^\#_1\cdot
i^\#_0)^\#, i^\#_1\cdot e(i^\#_0\cdot i^\#_1)^\#, i^\#_0\cdot
e(i^\#_1\cdot i^\#_0)^\#].
\]

 A functor between $B^\#$-categories is
called a {\it $B^\#$-functor} if it preserves the functors
$\times, \coprod$ and $\#$ and the fourteen natural
transformations $\Delta,p_0, p_1,
\nabla,i_0,i_1,\text{dist},m,j,e,l,r,\text{par}$ and $ \text{wd}$
defined above. Let $\mathbf{C}$ be a $B^\#$-category, and let $S$
be a set of morphisms of $\mathbf{C}$. The {\it $B^\#$-subcategory
generated by $S$} is the smallest $B^\#$-category of $\mathbf{C}$
containing $S$; this subcategory is denoted by $B^\#(S)$. If
$\mathbf{C}$ is a small $B^\#$-category, then the set
$\mathbf{C}_0$ of objects of $\mathbf{C}$ is an algebra with
signature $(\times,\coprod,\#)$; such an algebra is called a {\it
$(\times,\coprod,\#)$-algebra}. In the sequel, such algebras will
be obtained from certain objects of $B^\#$-category called
isotypical objects, to be defined. In particular, such an algebra
will generate a $B^\#$-category, in the following sense.

\begin{proposition}[~\cite{MR92e:03067}] \label{PropFreeBSharp}
If $\mathbf{C}$ is a small $B^\#$-category, $\mathbf{D}_0$ is a
$(\times,\coprod,\#)$-algebra and
$F_0:\mathbf{D}_0\rightarrow\mathbf{C}_0$ is a homomorphism, then
there exists, uniquely to within canonical isomorphism, a
$B^\#$-category $\mathbf{D}$ with object algebra $\mathbf{D}_0$,
supplied with a full and faithful $B^\#$-functor
$F:\mathbf{D}\rightarrow\mathbf{C}$ extending $F_0$.
\end{proposition}
\begin{proof}
Following the procedure in ~\cite{MR92e:03067} for obtaining the
smallest $B^\#$ category containing a given set of morphisms, we
think of the elements of $\mathbf{D}_0$ as the objects of an as
yet unspecified $B^\#$-category $\mathbf{D}$, and construct the
set $\mathbf{D}_1$ of morphisms accordingly. For each
$A\in\mathbf{D}_0$, we adjoin an identity morphism $1_A$ to
$\mathbf{D}_1$, subject to the functoriality relations $1_{A\times
B}=1_A\times1_B$, and so on. Next, adjoin the values of the
fourteen natural transformations above with arguments in
$\mathbf{D}_0$ to $\mathbf{D}_1$, subject to the relations i)
through v) above and subject to the relations that hold in a
$B^+$-category . Finally, we close the set of morphisms that
results under $\times, \coprod,\#$ and composition. This produces
the category $\mathbf{D}$. The homomorphism $F_0$ is extended to
$\mathbf{D}_1$ in the only way possible, following the three step
construction of $\mathbf{D}_1$. Identities in $\mathbf{D}_1$ must
be preserved by $F$. The image in $\mathbf{C}$ under $F$ of a
value of one of the fourteen natural transformations is completely
determined by functoriality and by definition of $F_0$. For
example, the value of $F$ on $\nabla_A:A\coprod A\rightarrow A$ is
$\nabla_{F A} :F A \coprod F A \rightarrow F A$. Finally, we
require that $F$ commute with $\times,\coprod, \#$ and
composition; for example, $F(\nabla_X\times
1_{Y^\#})=F\nabla_X\times 1_{F(Y^\#)}= \nabla_{F X}\times 1_{(F
Y)^\#}$.
\end{proof}

Let $\mathbf{C}$ be a $B^\#$-category. An object $X$  of
$\mathbf{C}$ generates by Proposition \ref{PropFreeBSharp} a
$B^\#$-category denoted by $\mathbf{C}\langle X\rangle$, whose
object algebra is free on the generator $X$, together with a
$B^\#$-functor $F:\mathbf{C}\langle X\rangle\rightarrow
B^\#(1_X)$. In case the $(\times,\coprod,\#)$-algebra generated by
$X$, namely $B^\#(1_X)_0$, is free, then $F$ gives an
identification of $\mathbf{C}\langle X\rangle$ with the full
subcategory $B^\#(1_X)$ of $\mathbf{C}$ and, following Heller, we
use the notation
\begin{align}\label{CBracketX}
\mathbf{C}=\mathbf{C}\langle X\rangle
\end{align}
to say that $\mathbf{C}_0$ is freely generated by $X$ as a
$(\times,\coprod,\#)$-algebra.

The definition of an iteration category requires the notion of a
stable union, a categorical notion which needs no special
definition in $\mathbf{Set}_F$. Let $\mathbf{C}$ be a
$P$-category, and let $X$ be an object of $\mathbf{C}$, and let
$\varepsilon_n\in\mathrm{dom}(X)$ be a countable family of
domains. A domain $\delta$ is a {\it union} of $\{\varepsilon_n\}$
if for each $n$, $\varepsilon_n \prec\delta$ and if, for any
$\phi$ and $\phi'$, $\phi\varepsilon_n=\phi'\varepsilon_n$ implies
$\phi\delta=\phi'\delta$. This determines $\delta$ uniquely, for
if $\delta'$ is any other union of $\{\varepsilon_n\}$, then
$\delta\varepsilon_n=\varepsilon_n=\delta'\varepsilon_n$, so
$\delta=\delta\delta=\delta'\delta$ and also
$\delta'=\delta'\delta'=\delta\delta'$. Since composition of
domains is commutative, $\delta=\delta'$. In this case. we set
$\delta=\cup\varepsilon_n$. If for any object $Y$ of $\mathbf{C}$,
\[
Y\boxtimes\bigcup\varepsilon_n=
\bigcup\left(Y\boxtimes\varepsilon_n\right),
\]
the union is called {\it stable} ~\cite{MR92e:03067}. In a
$P\Sigma$-category with ranges, every countable family of domains
$\{\varepsilon_n\}$ has the stable union
$\mathrm{ran}\left(\nabla\Sigma_{\mathbb{N}}\varepsilon_n\right)$.

An {\it iteration category} is a prodominical $P^+$-category with
ranges that satisfies the weak axiom of choice, such that the
$B^+$-category structure of $\mathbf{C}_T$ has been extended to a
$B^\#$-category so that for each object $X$, the formally free
semigroup $X^\#$ is the stable union of $\{\mathrm{ran}j_n\}$.

\begin{proposition}
$\mathbf{Pfn}_F$ is an iteration category.
\end{proposition}
\begin{proof}
Since $\mathbf{Pfn}_F$ is a $P\Sigma$-prodominical category with
ranges that satisfies the weak axiom of choice, it is an iteration
category.
\end{proof}

We give a translation of the iteration lemma from
Heller~\cite{MR92e:03067} into $\mathbf{Pfn}_F$. The iteration
lemma and its application to Turing data in the sequel is an
alternative to the use of fixed point semantics for iteration, as
exemplified by Manes' application of Kleene's fixed point theorem
to iterative specifications in $\omega$-complete categories
~\cite{MR94c:68134}. The proof of the iteration lemma that we give
is the transparent one available when the endofuctor $F$ preserves
products; in the general coalgebraic case we rely on the proof
given in Heller, which makes full use of the structure of a
$B^\#$-category, in lieu of our explicit manipulation of finite
sequences below.

\begin{lemma}[Iteration Lemma
~\cite{MR92e:03067}]\label{LemIteration} Let $X$ be a coalgebra,
let $i:U\hookrightarrow X$ be the inclusion of a sub-coalgebra $U$
into $X$, and let $f:X\rightarrow X$ be a morphism in
$\mathbf{Set}_F$ such that $f\circ i = i$. Then there is a unique
morphism $\mathrm{It}(f,U):X\rightarrow X$ in $\mathbf{Pfn}_F$
with domain
$\cup_n f^{-n}[U]$ such that for each %
$n\ge 0$, $\mathrm{It}(f,U)[f^{-n}[U]]=f^{n}[f^{-n}[U]]$.
Moreover, $\mathrm{im}\ \mathrm{It}(f,U)\subseteq U$ and, for
appropriate $g$ and $V$,
\begin{eqnarray*}
\mathrm{It}(f\times g,U\times V)=\mathrm{It}(f,U)\boxtimes\mathrm{It}(g,V),\\
\mathrm{It}(f\coprod g,U\coprod V) =
\mathrm{It}(f,U)\coprod\mathrm{It}(g,V).
\end{eqnarray*}
\end{lemma}
\begin{proof}
The proof we give in $\mathbf{Pfn}_F$ presupposes that the
endofunctor $F$ preserves products, since we assume that products
are constructed as they are in $\mathbf{Set}$; the general case
follows from ~\cite{MR92e:03067}.

Let $X^\#=\Sigma_\mathbb{N}X^{n+1}$ be the free semigroup
generated by $X$ and define the set
\[
W=\{(x,f x,\ldots,f^n x)\in X^\# : n\ge 1, x\in X, f^n x\in U\}.
\]
The set $W$ is a sub-coalgebra of $X^\#$ since, translating from
Heller ~\cite{MR92e:03067} into this context,
\[
W = \Delta_{X^\#}^{-1}\left[\left( (j\circ\mathrm{first}\cdot
f^\#)\times(X^\#\cdot
j\circ\mathrm{last})\right)^{-1}\left[\mathrm{im}\Delta_{X^\#}
\right]\right]\cap \mathrm{last}^{-1}[U],
\]
where $\mathrm{first}=[X,\ p_0]l_X$, $\mathrm{last}=[p_1,\ X]r_X$,
and where the dot $\cdot$ denotes concatenation of sequences in
$X^\#$; i.e., $(x_1,\ldots,x_m)\cdot(y_1,\ldots,y_n)=
(x_1,\ldots,x_m,y_1,\ldots,x_n)$.

Next we apply the diagram lemma in $\mathbf{Pfn}_F$. There exists
a unique induced map in the following diagram
\[
\xymatrix{ X^\# \ar[r]^-{\mathrm{first}|W} %
\ar[dr]_-{\mathrm{last}|W} & X
\ar@{-->}[d]^-{\frac{\mathrm{last}|W}{\mathrm{first}|W}} %
\\
& X }
\]
We may suppose that $n > m$. If
\[
\mathrm{first}(x,f x,\ldots,f^m x)= \mathrm{first}(y,f
y,\ldots,f^n y),
\]
then $x=y$, and as both of the indicated
sequences are in $W$, both $f^m x$ and  $f^n x$ are in $U$, and
therefore $f^{n-m}\circ f^m x = f^m x$. We take
\[
\mathrm{It}(f, U)=\frac{\mathrm{last}|W}{\mathrm{first}|W}.
\]
The remaining assertions are immediate.
\end{proof}

\subsection{Turing developments and local connectedness}
The iteration lemma \ref{LemIteration} will be applied to maps
$f:X\rightarrow X$ that fix a sub-coalgebra $U$ arising as a
summand of a coproduct $X = U\coprod V$.

Let $\mathbf{C}$ be an iteration category. A {\it Turing datum} in
$\mathbf{C}$ is a diagram
\[
\xymatrix{X\ar[r]^-u & W\ar[r]^-v\ar[r]& W\coprod Y }
\]
in $\mathbf{C}_T$; in our case  $\mathbf{C}=\mathbf{Pfn}_F$
and  $\mathbf{C}_T=\mathbf{Set}_F$. The map %
$[v, i_1]:W\coprod Y\rightarrow W\coprod Y$ is total and satisfies %
$[v, i_1]\left(\emptyset\coprod Y\right)=\emptyset\coprod Y$, so
by the iteration
lemma \ref{LemIteration} there is a map %
\[
\xymatrix{ \mathrm{It}\left([v, i_1],\emptyset\coprod
Y\right):W\coprod Y \ar[r] & W\coprod Y}.
\]
 The {\it Turing development}
$\mathrm{Tur}(u,v):X\rightarrow Y$ of the given Turing datum is
the composite
\[
\xymatrix{X\ar[r]^-{i_0 u} & W\coprod Y\ar[rrr]^-{\mathrm{It}([v, i_1], \emptyset\coprod Y)}%
&&& W\coprod Y \ar[r]^-{[0, 1_Y]}& Y}
\]
If $\mathbf{D}$ is a $B^+$-subcategory of $\mathbf{C}_T$, then the
{\it class of all Turing developments} of Turing data in
$\mathbf{D}$ is denoted by $\mathrm{Tur}(\mathbf{D})$. It follows
from Lemma 8.1 of Heller ~\cite{MR92e:03067} that
$\mathbf{D}\subseteq \mathrm{Tur}(\mathbf{D})$ and that
$\mathrm{Tur}(\mathbf{D})$ is closed under $\boxtimes$ and
$\coprod$. Under the additional assumption that $\mathbf{C}$ is
``locally connected", it follows that $\mathrm{Tur}(\mathbf{D})$
is closed under composition ~\cite{MR92e:03067}.

A coalgebra is {\it connected} if it is connected in the coalgebra
topology; i.e., if it is not a nontrivial coproduct of
subcoalgebras. If for each coalgebra $X$ of the iteration category
$\mathbf{Pfn}_F$, the coalgebra topology on $X$ is locally
connected, we say that $\mathbf{Pfn}_F$ is {\it locally
connected}; in that case by Lemma 8.2 of Heller
~\cite{MR92e:03067}, for any $B^+$ subcategory $\mathbf{D}$ of
$\mathbf{C}$, $\mathrm{Tur}(\mathbf{D})$ is closed under
composition; moreover, it is a $+$-prodominical subcategory of
$\mathbf{C}$ containing $\mathbf{D}$.

\begin{proposition}
If the $\mathbf{Set}$ endofunctor $F$ weakly preserves generalized
pullbacks of monomorphisms, then each $F$-coalgebra is locally
connected in the coalgebra topology.
\end{proposition}
\begin{proof}
Under the hypothesis on $F$ it follows from Theorem 5.10 of Gumm
and Schr\"oder ~\cite{MR1787576} that for each  $F$-coalgebra $X$
and for each $x\in X$, the subcoalgebra cogenerated by $x$ exists;
hence the $1$-cogenerated subcoalgebras of $X$ form a base for the
coalgebra topology on $X$. Each $1$-cogenerated subcoalgebra is
connected in the coalgebra topology.
\end{proof}

\subsection{Isotypical objects and isotypes}
We paraphrase Heller ~\cite{MR92e:03067}. A category $\mathbf{C}$
is {\it isotypical} if any two of its objects are isomorphic. For
example, if $\mathfrak{M}$ is an infinite cardinal number, then
the full subcategory $\mathbf{Set}_\mathfrak{M}$ of $\mathbf{Set}$
containing the sets of cardinality $\mathfrak{M}$ can be given the
structure of a $B\Sigma$-category, hence of an isotypical
$B^\#$-category.

Let $\mathbf{C}$ be a $B^\#$ category. An object $X$ of
$\mathbf{C}$ is {\it isotypical} if it is isomorphic to each of
$X\times X$, $X\coprod X$, and $X^\#$. Such objects are used in
the construction of recursion categories in the sequel. If $X$ is
an isotypical object of a $B^\#$-category $\mathbf{C}$, then
$\mathbf{C}\langle X\rangle$ (cf. equation \ref{CBracketX}) is an
isotypical $B^\#$-category.

\subsection{Uniform generation}
The notion of the domain of a morphism carries over from
$P$-categories to $B$-categories; in particular, if $\mathbf{C}$
is a $B$-category, then the notion of a total morphism is
definable in $\mathbf{C}$. If $\mathbf{C}$ is a $B$-category and
$t:W\times X\rightarrow X$ is a morphism in $\mathbf{C}$, then an
{\it index} of $f:X\rightarrow X$ relative to the {\it catalog}
$t$ is a total morphism $g:X\rightarrow W$ such that the following
diagram commutes.
\begin{eqnarray}\label{GeneralizedIndex}
\xymatrix{X\times X \ar[d]_{p_1} \ar[r]^-{g\times 1} & W\times X
\ar[d]^-t\\
X\ar[r]^-f & X}
\end{eqnarray}
If $\mathbf{B}$ is a subcategory of $\mathbf{C}$, the {\it uniform
list cataloged by $t$ with indices in $\mathbf{B}$} is the set
$\mathcal{L}(\mathbf{B},t)$ of morphisms $f\in\mathbf{C}(X,X)$ for
which there exists a total $g\in\mathbf{B}(X,W)$ such that the
diagram (\ref{GeneralizedIndex}) commutes.

If $\mathbf{C}$ is a $B^\#$-category and
$\mathbf{B}=\mathbf{B}\langle X\rangle$ is an isotypical
$B^\#$-subcategory, then $\mathbf{B}$ contains a {\it frame $b$ at
$X$}, namely a collection of isomorphisms $b_\times:X\rightarrow
X\times X,b_{\coprod}:X\rightarrow X\coprod X,b_\#:X\rightarrow
X^\#$, along with their inverses. The subcategory $\mathbf{B}$ is
called a {\it uniformly generated} $B^\#$-subcategory of
$\mathbf{C}\langle X\rangle$ provided
\begin{align}
\mathbf{B}=B^\#\left(b\cup \mathcal{L}(\mathbf{B},t)\right)
\label{Uniform}
\end{align}
(cf. Proposition \ref{PropFreeBSharp} and preceding remarks).

Uniformly generated isotypical categories have the following
properties. If $b$ and $b'$ are frames at $X$ contained in the
uniformly generated category $\mathbf{B}$, then
\[
\mathbf{B}=B^\#\left(b\cup \mathcal{L}(\mathbf{B},t)\right)=
B^\#\left(b'\cup \mathcal{L}(\mathbf{B},t)\right).
\]
Any finitely generated isotypical $B^\#$-subcategory
$\mathbf{B}=\mathbf{B}\langle X\rangle\subset\mathbf{C}\langle
X\rangle$ is uniformly generated. For any $t:W\times X\rightarrow
X$ in $\mathbf{C}\langle X\rangle$ there is a uniformly generated
$\mathbf{B}\subset\mathbf{C}$ containing $t$ and
$\mathcal{L}(\mathbf{B},t)$. For any frame $b$ and for any
catalogue $t$ such that (\ref{Uniform}) holds, there exists a
maximal uniformly generated category $\mathbf{B}$ satisfying
(\ref{Uniform}).

\section{Recursion categories}

\subsection{The existence theorem}
A {\it Turing morphism} in a prodominical isotype is a morphism %
$\tau:W\boxtimes X\rightarrow Y$ such that for any
$\phi:V\boxtimes X\rightarrow Y$ there exists a total
$g:V\rightarrow W$ such that the following diagram commutes.
\[
\xymatrix{V\boxtimes X \ar[r]^-{g\times 1_X} \ar[dr]_-\phi &
W\boxtimes X
\ar[d]^-\tau \\
& Y }
\]
A {\it recursion category} is a prodominical isotype in which
there is a Turing morphism.

The statement of Heller's existence theorem follows.
\begin{theorem}[Heller ~\cite{MR92e:03067}]
Let $\mathbf{C}$ be a locally connected iteration category. If
$\mathbf{D}=\mathbf{D}\langle X \rangle$ is a uniformly generated
isotypical $B^\#$-subcategory of $\mathbf{C}\langle X \rangle$,
then $\mathrm{Tur}(\mathbf{D})$ is a recursion category.
\end{theorem}

By preceding remarks, we have the following.
\begin{theorem}
Let $F$ be a nontrivial endofunctor on $\mathbf{Set}$ that weakly
preserves pullbacks. Suppose that $\mathbf{Set}_F$ is locally
connected and complete. Then for any uniformly generated
isotypical subcategory $\mathbf{C}$ of $\mathbf{Pfn}_F$,
$\mathrm{Tur}(\mathbf{C})$ is a recursion category.
\end{theorem}

\begin{corollary}
Let $F$ be a nontrivial bounded or accessible endofunctor on
$\mathbf{Set}$ that weakly preserves generalized pullbacks. If $X$
is an isotypical coalgebra of $\mathbf{Pfn}_F$ and if $b$ is a
frame at $X$, then for any finite collection of morphisms
$S\subset\mathbf{Pfn}_F(X,X)$, $\mathrm{Tur}\left(B^\#(b\cup
S)\right)$ is a recursion category.
\end{corollary}

\section{Examples and remarks}

Taking the $\mathbf{Set}$ endofunctor $F$ to be the identity
functor and each coalgebra structure map $\alpha:X\rightarrow F X$
to be the identity morphism $1_X$, one obtains the category
$\mathbf{Pfn}$ of sets and partial functions as a subcategory of
$\mathbf{Pfn}_F$; hence $\mathbf{Pfn}_F$ yields all the examples
of recursion categories that come from sets and partial functions,
as in ~\cite{MR92e:03067}. The set of functions computable in the
Blum-Shub-Smale model of computation ~\cite{MR99a:68070} over a
ring can be obtained as morphisms of an appropriately defined
recursion category. We omit the construction.

Since we suppose that the endofunctor $F$ is such that
$\mathbf{Set}_F$ is complete, e.g. if $F$ is bounded or
$\omega$-accessible, $\mathbf{Set}_F$ contains a terminal object
$1$. For any $F$-coalgebra $X$, let $Y=\coprod_{\mathbb{N}}X$ be a
countable copower, let $Z=\coprod_{\mathbb{N}}Y^n=1\coprod
Y\coprod Y^2\coprod\cdots$ and let $W=Z^\#$. By proposition 1.1 of
~\cite{MR92e:03067}, $W$ is an isotypical object of
$\mathbf{Set}_F$ and hence of $\mathbf{Pfn}_F$. If, in addition,
$F$ preserves weak limits of sinks, then $\mathbf{Pfn}_F$ is
locally connected, and by the theorem its isotypical subcategory
$\mathbf{Pfn}_F\langle W\rangle$ contains many recursion
categories.

We remark on further aspects of $P$-categories of coalgebras
beyond the scope of this paper. On the assumption that $F$ be a
nontrivial endofunctor on $\mathbf{Set}$ that weakly preserves
pullbacks and such that $\mathbf{Set}_F$ is complete,
$\mathbf{Pfn}_F$ is a ranged Boolean category. Boolean categories
were defined by Manes ~\cite{MR94c:68134} as a categorical setting
for predicate transformer semantics. Predicate transformers were
introduced by Dijkstra and have been applied as a formal calculus
of program derivation from specifications in logic
~\cite{MR91c:68089,MR93a:68001}. It follows from Manes' theory
that $\mathbf{Pfn}_F$ admits a representation as a category of
relations.

\subsection*{Acknowledgement}
I wish to thank Alex Heller for introducing me to recursion
categories and for helpful remarks, Noson Yanofsky for carefully
reading a draft of this paper and for suggestions, and Ellis
Cooper for discussions on coalgebras and for editorial comments.

\bibliographystyle{amsalpha}
\bibliography{rec30May01}
\end{document}